%% file: paper_revised.tex
\newcommand{\C}{\mathbb{C}}
\newcommand{\Z}{\mathbb{Z}}
\input{defns}
\documentclass[reqno]{amsart}
\usepackage{amssymb, amsthm, amsmath, latexsym}
\usepackage{epsfig}

\title{Minimal Tori in $S^3$.}
\author{Emma Carberry}
\address{Department of Mathematics\\University of Sydney
 \\ NSW, 2006, Australia} \email{carberry@maths.usyd.edu.au}


\theoremstyle{plain}
\newtheorem{thm}{Theorem}[section]
\newtheorem{lem}[thm]{Lemma}

\newtheorem{corollary}[thm]{Corollary}
\newtheorem {defn}[thm]{Definition}

\newtheorem* {thm*}{Theorem}
\newtheorem *{cor*}{Corollary}
\newtheorem{conditions}[thm]{Periodicity Conditions}

\begin{document}

\begin{abstract}
We prove existence results that give information about the space of minimal immersions of 2-tori into $ S ^ 3 $. More specifically, we show that
\begin{enumerate}
\item For every positive integer $ n $, there are countably many real $n $-dimensional families of minimally immersed 2-tori in $ S ^ 3 $.  Every linearly full minimal immersion $ T ^ 2\rightarrow S ^ 3 $ belongs to exactly one of these families.
\item Let $ \mathcal A $ be  the space of rectangular 2-tori. There is a countable dense subset $\mathcal B  $ of $\mathcal A  $  such that every torus in $\mathcal B$ can be minimally immersed into $ S ^ 3 $.
\end{enumerate}
The main content of this manuscript lies in finding minimal immersions that satisfy {\bf periodicity conditions} and hence obtaining maps of tori, rather than simply immersions of the plane.
We make use of a correspondence, established by Hitchin, between minimal tori in $S^3$ and algebraic curve data.
\end{abstract}
\maketitle
\tableofcontents
\section{Introduction}
We show that there is an abundance of minimally immersed 2-tori in $S^3$; in fact they come in families of every real dimension.
 We utilise results of Hitchin, \cite{Hitchin:90}, namely his description of harmonic maps $f:T^2\to S^3\simeq
SU(2)$ in terms of algebraic curve data. This transforms a problem in analysis to one in algebraic geometry, hence rendering it vulnerable to new techniques.  A conformal map is harmonic if and only if it is minimal, and so this reformulation applies to minimal immersions.  The  challenge of finding {\bf doubly-periodic} minimal immersions (ie immersions of tori) hence becomes a problem in algebraic geometry, and it is this problem that we solve.  Our approach builds upon work of Ercolani, Kn{\"o}rrer and Trubowitz \cite{EKT:93}, in which they proved it is possible to obtain families of constant mean curvature tori in $\R ^ 3 $ of every even dimension. Hitchin's correspondence excludes minimal branched immersions which are not linearly full, i.e. which map into a totally geodesic $S^2\subset S ^ 3$.  We henceforth assume that our minimal immersions are linearly full.  We also partially address the interesting question of which tori can be minimally immersed into the 3-sphere. In particular we show that a countable dense subset $\mathcal B$ of the space of rectangular tori do allow such immersions. 

Exploiting the algebro-geometric viewpoint allows us to obtain not only strong existence results, but also  detailed information about the {\bf space $\mathcal M $
of linearly full minimal immersions of $ T ^ 2\rightarrow S^3\cong SU(2)$}. These are one of the most basic examples of harmonic maps from a torus to a Lie group or symmetric space, and provide a model for the general situation.  One can also associate a spectral curve to harmonic maps of tori to other  Lie groups and symmetric spaces, although the general theory is incomplete  (see eg \cite{Griffiths:85, Burstall:92, BFPP:93, FPPS:92, McIntosh:95, McIntosh:96}). There has been an explosion of interest in these maps over the last 30 years, both due to their geometric interest and their strong connections with the Yang-Mills equations. These results can be viewed as a preliminary step in the very interesting general program of obtaining information about the spaces of these harmonic maps.

%
We now explain how the algebro-geometric approach elucidates the the structure of  $\mathcal M $, and why it is a useful approach to the study of these (and other) maps.
In \cite {Hitchin:90},  Hitchin established an explicit bijective correspondence between harmonic maps  $ T ^ 2\rightarrow S ^ 3 $ and {\bf
spectral curve data}, which includes a hyperelliptic curve $\Sigma$, called a spectral curve, and a line bundle $ E (0)\rightarrow\Sigma $.  (See Definition~\ref{spectraldata} and Theorem~\ref{thm:Hitchin} for details). This introduces an important new invariant of the harmonic map, namely the arithmetic genus $ g $ of $\Sigma $, or {\bf spectral genus}.
The space $\mathcal M $ decomposes into strata $\mathcal S_g$, where $\mathcal S _ g $ consists of  maps of spectral genus $ g $. We show that for $ g> 2 $, $\mathcal S _ g $ is comprised of (at least) countably many  $ (g - 2) $-dimensional  families of immersions. From \cite{Hitchin:90}, $ \mathcal S _ 1 $ is the $ S ^ 1 $-symmetric examples, and $\mathcal S _ 0 $ is just the finite covers of the Clifford torus.

It is easy to write down minimal immersions of the plane to $S^3$, but finding ones which are doubly periodic is a difficult problem.  Hitchin's correspondence transforms this  analytic problem into a similarly non-trivial algebro-geometric one.  Namely, it becomes the requirement that $\Sigma$ supports a pair $\Theta$, $\Psi$ of meromorphic differentials of the second kind, whose periods,  together with their integrals over
certain open curves, $C_1$ and $C_{-1}$, are all  integers. These {\bf periodicity conditions} place very strong restrictions on the spectral data; they demand that certain transcendental functions are integer valued, and a priori it is not at all clear that solutions exist for each $ g $.

Hitchin proved the existence of harmonic maps $T^2\to S^3$ of spectral genera $g\leq 3$, and of minimal tori with $g\leq 2$. Our main result is Theorem~\ref{thm:main},
which gives the algebro-geometric statement necessary for our  conclusions regarding the space $\mathcal M$ of minimally immersed tori. It states that
\begin{thm*}[\bf{\ref{thm:main}}]
For each integer $g>0$ there are countably many spectral curves of arithmetic genus $g$ that give rise to minimal
immersions from rectangular tori to $S^3$.
\end{thm*}
Utilising Hitchin's spectral curve correspondence, we conclude
\begin {cor*}[\bf{\ref{thm:families}}] For each integer $n\geq 0 $ there are countably many real $ n $-dimensional families of minimal immersions from rectangular tori to $ S ^ 3 $.  Each family consists of maps from a fixed torus.
\end {cor*}
Bryant has proven that all compact Riemann surfaces can be minimally immersed into $ S ^ 4 $ \cite {Bryant:82}, but  the question of which Riemann surfaces or indeed which tori admit minimal immersions into $ S ^ 3 $ is both open and very interesting. The area of a surface is invariant under conformal transformations, and so one only needs a conformal structure on a surface to have a notion of minimal immersions of it.  Thus when we ask `` which tori?" we mean " which conformal classes of tori?".

Using conformal equivalence and taking covers, the following corollary is straightforward.
\begin {cor*}[\bf{\ref{thm:conftype}}]
Let $ \mathcal A$ denote  the space of rectangular 2-tori,  and call a torus {\it admissible} if it admits a linearly full minimal immersion into $ S ^ 3 $.  Write $\mathcal C$ for the space of admissible rectangular 2-tori, and $\mathcal C _ g $ for those which possess a minimal immersion of spectral genus $g $. Then for each $ g\geq 0$, $\mathcal C_g$ is (at least) countable and dense in $\mathcal A $.
%
\end {cor*}


We now briefly discuss the methods by which we prove Theorem~\ref{thm:main}. We impose additional symmetry conditions on our spectral data in order to ensure that the resulting tori are rectangular.  Using these and other symmetries, one can define a smooth map $\phi $ from an open subset of the space of spectral curves of genus $ g $ to a  space of the same dimension in which the periodicity conditions are satisfied on a countable dense subset. The technical heart of the problem  is then to show that for each $ g $, there is a curve of genus $ g $ at which the differential of this map $\phi $ is invertible.

%

Our methods  extend those of Ercolani, Kn{\"o}rrer and Trubowitz
\cite{EKT:93}, who proved the existence of rectangular constant
mean curvature tori in $\R^3$ for every even spectral genus $g\geq
0$ (there is also a spectral curve description of such tori). The
additional symmetry that forces the resulting tori to be
rectangular results in the proof breaking naturally into odd and
even genus cases. Our main contribution is the odd genus case,
which is more complicated and is presented in detail.  For even
genera the proof is reminiscent of \cite{EKT:93} and is briefly
indicated. We remark that Jaggy   \cite{Jaggy:94} showed the
existence of constant mean curvature tori in $\R ^ 3 $ for all
spectral genera, but we do not pursue analogous methods since they
lead to less refined information concerning the conformal type. It
is by imposing this additional symmetry that we are able to obtain
the conformal type information stated in
Corollary~\ref{thm:conftype}.

We now mention some open questions that are motivated by this manuscript.
\begin {enumerate}
\item We prove that a countable dense subset of the space of rectangular tori can be minimally immersed in $ S ^ 3 $.  Our method of proof shows that  in an open subset of the space  $\mathcal M $
of linearly full minimal immersions of $ T ^ 2\rightarrow S^3$ these are  the only rectangular tori that occur, but this result is only local. The spectral curve viewpoint suggests that the space $\mathcal C$ of admissible rectangular tori is only countable, and hence that in contrast to the situation for $ S ^ 4 $, not all tori can be minimally immersed into $ S ^ 3 $. One could prove this by showing that the differential of $\phi $ is everywhere invertible.
\item The question of exactly which tori lie in $\mathcal C  $ is  difficult, but very interesting. Further, are the $ \mathcal C _ g $ different for different $ g $? It is not even known whether the square torus lies in all $\mathcal C _ g $.
\item These questions have natural analogues for general (not necessarily rectangular) tori.
\item An obvious problem that this work suggests is to prove similar existence results for ``non-rectangular" tori. Some care needs to be exercised in stating this problem, since we can always replace a given torus by a conformally equivalent one, what we are seeking here are tori which are not conformally equivalent to any rectangular torus.  Methods similar to those employed here will yield non-rectangular tori (cf  \cite{Jaggy:94} ).  However, it is not clear whether these tori will  be conformally equivalent to rectangular tori; the task is to show that for each $ g $, there are some which are not.
\end {enumerate}
We conclude the introduction by outlining the contents of this paper.  In section 2, we provide a short account of Hitchin's spectral curve correspondence and give our geometric reformulation of the periodicity conditions. Next, in section 3, we state our main results and explain the strategy for their proof. Section 4 is the technical heart of this manuscript, and is where the proof is carried out.

The results of this paper are taken from the author's thesis
\cite{Carberry:02}, in which further details can be found. It is a pleasure to acknowledge here the invaluable support and wisdom of my advisor, Phillip A. Griffiths, as
well as insightful conversations with Chuu-Lian Terng, Ivan Sterling and Ian McIntosh.

\section{Background: Spectral Curve Correspondence}
In this section, we first explain how a harmonic map from a surface to a Lie group with bi-invariant metric can be described in terms of a family of flat connections.  We then outline Hitchin's spectral curve construction for harmonic maps from a 2-torus to $ S ^ 3\cong SU(2)$, which further reduces the harmonic map equations to a linear flow in a complex torus, namely the Jacobian of the spectral curve.

Let $M$ be a compact Riemann surface, $G$ a compact Lie group with bi-invariant metric, and
$\mathfrak{g}$ its Lie algebra. Writing $\Phi$ for the $(1,0)$ part of $f^{-1}df$,
\[ f ^ {- 1}df = \Phi - \Phi^*.\]
Here $\Phi$ is valued in the complex Lie algebra with real form $\mathfrak{g}$, and $-\Phi^*$ is its image under
the corresponding real involution. For $G$ a unitary group, $\mbox {$\Phi^*=\bar\Phi^t$} $.

 A smooth map $f: M\rightarrow G$ is harmonic if and only if
%
\begin{equation} d\phi_\lambda+\frac 12[\phi_\lambda,\phi_\lambda]=0 \mbox{ for all }
\lambda\in\C^*,\label{eq:flat}\end{equation}
where
\[
\phi_\lambda:=\frac 12 (1-\lambda^{-1})\Phi - \frac 12
(1-\lambda)\Phi^* ,\,
\lambda\in\C^*.
\]
Let $P$ be a trivial principal $SU(2)$-bundle over $M$, with trivial connection
$\nabla$.  Writing
\begin {equation}\nabla_\lambda:=\nabla + \phi_\lambda, \label{eq:conn}\end {equation}
$f$ gauges the trivial connection \mbox{$\nabla_{1}=\nabla$} to $\nabla_{-1}$, and (\ref{eq:flat}) states that
$\nabla_\lambda$ is flat for each $\lambda\in\C^*$.

Conversely, take a principal $SU (2) $-bundle $P$ , a
fixed connection $\nabla$   and $\Phi\in\Omega^{1,0}(M,Ad(P)\otimes\C)$ such that the connections $\nabla_\lambda$
defined as above are flat. Then locally one recovers a harmonic map $f:U\to G$ by
\begin {equation} g_{-1}(z)=f(z)g_{1}(z),\label {eq:sections}\end {equation}
where $g_1, g_{-1}$ are sections of $P$ over $U$ parallel with
respect to $\nabla_1$ and $\nabla_{-1}$ respectively; patching
gives a harmonic section of $P$. The condition that we obtain a global $f:M\to G$  is
that both the left connection $\nabla_{1}=\nabla$ and the right
connection $\nabla_{-1}$ are trivial, and then $f$ is determined up
to left and right actions of $G$.

We outline Hitchin's algebro-geometric description of harmonic maps $f:T^2\to SU(2)$,  \cite{Hitchin:90}. We can consider the
flat connections $\nabla_\lambda$ as living in a rank two complex vector bundle $V$ with a symplectic form
$\omega$ and a quaternionic structure.  $f$ is determined from the connections $\nabla_\lambda$ up to the action of $SO(4)=SU(2)\times SU(2)$
on $S^3$.  The restriction to tori is essential, as the construction
requires the compact Riemann surface to have abelian (and non-trivial) fundamental group.
It is
assumed throughout that $f$ is not a conformal harmonic map to a totally geodesic \mbox{$S^2\subset S^3$}; such
maps do not admit a spectral curve. First, one gives an algebro-geometric description of families of flat connections (\ref {eq:conn}), or equivalently to harmonic sections $ f $ of $ P $.

Given a family (\ref{eq:conn}) of flat
connections on a marked torus $T^2$, let $H_z(\lambda), K_z(\lambda) \in SL(2,\C)$ denote the holonomy of
$\nabla_\lambda$ for the chosen basis \mbox{$[z,z+1]$}, \mbox{$[z,z+\tau]$} of \mbox{$\pi_1(T^2, z)$}, and let
$\mu(\lambda)$, $\nu(\lambda)$ be their eigenvalues.

Define the spectral curve $\Sigma$:
\[
\eta^2=p(\lambda)=\prod_{i=1} ^ {n} (\lambda -\lambda _ i)
\]
as follows.  For generic $\lambda\in\C^*$,  $ H _ z (\lambda)$ is not $\pm I$, and so has two eigenvectors $ v _ 1(\lambda), v _ 2 (\lambda) $, which may coincide.  $\lambda $ is a zero of $p$ of order $k$ if $ v _ 1(\lambda), v _ 2 (\lambda)\in V _ z $ agree to order $ k $, as measured by the order of  vanishing of the symplectic form $\omega $.   If $ H _ z (\lambda)=\pm I$, one can use holomorphic continuation to choose $ v _ 1(\lambda) $ and $ v _ 2 (\lambda) $, and then apply the same criterion. See \cite {Hitchin:90} for a more careful treatment.
If $f$ is conformal (det $\Phi $ = 0) then $p$ has odd order and
has a simple zero at $0$ whilst if $f$ is nonconformal then $p$ has even order and is non-zero at $0$.

Since the eigenspaces of the holonomy matrices
$H_z(\lambda)$ for different $z$ are related by parallel transport, $\Sigma $ is independent of $z$. It follows from the abelianicity of $\pi _ 1 (T ^ 2) $ that it is also independent of the choice of generator of the fundamental group, and crucially, it is a (finite genus) algebraic curve \cite {Hitchin:90}.

Clearly the eigenvalue functions $\mu$ and $\nu$ are well-defined regular
functions on $\Sigma -\lambda^ {- 1} (0,\infty) $. $$\Theta:=\frac{1}{2\pi i}d\log\mu\mbox{,
}\Psi:=\frac{1}{2\pi i}d\log\nu$$ are differentials of the second kind on $\Sigma$; their only singularities are double poles at
$\pi^{-1}\{0,\infty\}$, they have no residues and their periods are integers.

For each $z\in T^2$, ${\Sigma}-\pi^{-1}\{0,\infty\}$ supports the
eigenspace line bundle
\[\left(E(z)\right)_{(\lambda,\eta)}\subseteq
\ker(H(\lambda,z)-\mu(\lambda,\eta))\subseteq V_z,\] which extends to a holomorphic line bundle $E(z)$ on
${\Sigma}$, and the map $z\rightarrow E(z)$ is linear. $\Sigma $ has a hyperelliptic involution $\sigma $, and the $SU(2)$ structure also induces an antiholomorphic involution $\rho $, covering $\lambda\mapsto\bar{\lambda}^{-1}$.

Thus by explicit construction, Hitchin associates the following  spectral data
$(\Sigma,\lambda, \rho, \Theta, \Psi, E(0))$ to the family
(\ref{eq:conn}). %

\begin{defn} [ Spectral Data, \cite {Hitchin:90}, Theorem 8.1]\label{spectraldata}
By spectral data we mean that
\begin{enumerate}

\item $\Sigma$ is a hyperelliptic curve $\eta^{2}=p(\lambda)$ of arithmetic genus $g$, and
\mbox{$\pi:\Sigma\to\mathbb{C}P^{1}$} is the projection \mbox{$\pi(\lambda,\eta)=\lambda$.}

\item $p(\lambda)$ is real with respect to the real structure
\mbox{$\lambda\mapsto\bar{\lambda}^{-1}$} on $\mathbb{C}P^{1}$.

\item $p(\lambda)$ has no real zeros (i.e. no zeros on the unit
circle $\lambda = \bar\lambda^{-1}$).

\item $p(\lambda)$ has at most simple zeros at $\lambda=0$ and
$\lambda=\infty$.

\item $\Theta$ and $\Psi$ are meromorphic differentials on
$\Sigma$ whose only singularities are double poles at $\pi^{-1}(0)$ and
$\pi^{-1}(\infty)$ and which have no residues. Their principal parts are linearly
independent over $\mathbb{R}$, and they satisfy
\begin{center}$\sigma^{*}\Theta=-\Theta,\,\sigma^{*}\Psi=-\Psi,\,
\rho^{*}\Theta=\bar\Theta,\,\rho^{*}\Psi=\bar\Psi$\end{center}
where $\sigma$ is the hyperelliptic involution
$(\lambda,\eta)\mapsto(\lambda,-\eta)$ and $\rho$ is the real
structure induced from $\lambda\mapsto\bar\lambda^{-1}$.

\item {\it The periods of $\Theta$ and $\Psi$ are all  integers.}

\item $ E (0) $ is a line bundle of degree $g+1$ on $\Sigma$, quaternionic
with respect to the real structure $\sigma\rho$.

\end{enumerate}
\end {defn}

 Conversely, from spectral data $(\Sigma,\lambda, \rho, \Theta, \Psi,
E (0))$ one can recover a marked torus $T^2$, an $SU(2)$ principal bundle $P$, and a family of
 flat connections of the form (\ref{eq:conn})
 unique up to  gauge transformations and the operation of
tensoring the associated complex vector bundle $V$ by a flat $\mathbb{Z}_{2}$-bundle (Theorem 8.1,
\cite{Hitchin:90}).

The existence of differentials $\Theta$, $\Psi$ with integral periods places a stringent constraint on $\Sigma$.
It corresponds to the double periodicity of a harmonic section of $P={G\times G}\diagup{G}$.

The following is essentially Theorem 8.20 of \cite{Hitchin:90}. The additional periodicity conditions guarantee
that the connections $\nabla _ 1, \nabla _ {-1} $ are trivial and hence the harmonic section (\ref{eq:sections}) of $P$ gives a  map $ T ^ 2\rightarrow S ^ 3 $.
\begin{thm} [Hitchin, \cite{Hitchin:90}] \label{thm:Hitchin}
Let $(\Sigma,\lambda, \rho,\Theta,\Psi,E(0))$ be spectral data, $\mu$
and $\nu$ be functions on $\Sigma-\lambda^{-1}\{0,\infty\}$ satisfying
$\Theta=\frac{1}{2\pi i}d\log\mu,\,\Psi=\frac{1}{2\pi i}d\log\nu $
and $\mu\sigma^*\mu=1,\,\nu\sigma^*\nu=1$ (such functions exist by
the periodicity conditions on $\Theta$, $\Psi$). Let $T^2$ be the marked torus with generators \mbox{$[z,z+1]$}, \mbox{$[z,z+\tau]$} for $\pi_1(T^2,z)$ where $\tau=\frac{\mbox{principal part}_{P_0}\Psi}{\mbox{principal part}_{P_0}\Theta}$ (in some local coordinate system) and $P_0\in\lambda^{-1}(0)$. Then
\begin{enumerate}
\item $(\Sigma,\lambda, \rho,\mu,\nu, E(0))$ determines a harmonic map
$f:T^2\to SU(2)$ if and only if
$$\mu(\lambda,\eta)=\nu(\lambda,\eta)=1\mbox{ for all
}(\lambda,\eta)\in\lambda^{-1}\{1,-1\}.$$

\item $f$ is conformal if and only if $P(0)=0$.


\item The harmonic map is uniquely determined by $(\Sigma,\lambda,
\rho, \mu,\nu, E(0))$ modulo the action of $SO(4)$ on $S^3$.
\end{enumerate}
\end{thm}

Let $C_1$ be a curve in $\Sigma$ joining the two points in $\lambda^{-1}(1)$, and  $C_{-1}$  a curve joining
the two points  in $\lambda^{-1}(-1)$.
\begin {conditions}\label {periodicity}
The existence of functions $\mu$, $\nu$ as above is equivalent to the following periodicity conditions ($\mu$ and $\nu$ are then determined up to sign).
\begin{enumerate}
\item The periods of $\Theta,\Psi $ are all integers.
\item\[\int_{C_{\pm 1}}\Theta,\,\int_{C_{\pm 1}}\Psi\in\mathbb Z,\]
\end{enumerate}
\end {conditions}
%
%
%

\section{Statement of Results}

Here we state the main algebro-geometric result that we prove (\ref{thm:main}), and explain how it enables us to obtain information regarding the space of minimally immersed tori in $ S ^ 3 $ (Corollaries \ref{thm:families}, \ref{thm:conftype}). We also outline our approach to the proof of this theorem.

We will prove the following theorem:

\begin{thm}\label{thm:main}
For each integer $g>0$ there are countably many spectral curves of arithmetic genus $g$ that give rise to minimal
immersions from rectangular tori to $S^3$.
\end{thm}

%
\begin {corollary} \label{thm:families} For each integer $n\geq 0 $ there are countably many real $ n $-dimensional families of linearly full minimal immersions from rectangular tori to $ S ^ 3 $.  Each family consists of maps from a fixed torus.
\end {corollary}

\begin {proof} Given spectral data
$(\Sigma, \lambda, \rho, \Theta, \Psi)$, condition
\setcounter{enumi}{3}({\it\arabic{enumi}}\,) of Hitchin's spectral data (Definition~\ref {spectraldata}) implies that the quaternionic structure $\rho\sigma$ has no fixed points and so by
\cite{Atiyah:71} there is a real $g$-dimensional family of quaternionic bundles of degree $g+1$. Factoring out by reparameterisations of the domain torus, for $g>2$ we obtain a real  $(g-2)$-dimensional family of minimal immersions from a fixed torus.
\end {proof}

The following corollary is then  straightforward.
\begin {corollary}\label{thm:conftype}
Let $ \mathcal A$ denote  the space of rectangular 2-tori,  and call a torus {\it admissible} if it admits a linearly full minimal immersion into $ S ^ 3 $.  Write $\mathcal C$ for the space of admissible rectangular 2-tori, and $\mathcal C _ g $ for those which possess a minimal immersion of spectral genus $g $. Then for each $ g\geq 0$, $\mathcal C_g$ is (at least) countable and dense in $\mathcal A $.
 \end {corollary}

\begin {proof} $\mathcal C_0$ contains the square torus (Clifford torus) and we have shown above that each $\mathcal C _ g $, $ g>0 $ is non-empty. Let $\tau $ be such that the torus with sides $ [ 0, 1] $ and $ [0,\tau ] $ is in $\mathcal C _ g $.  Then  for each $ q\in\mathbb {Q} $, the torus with sides $ [ 0, 1] $ and $ [0,q\tau ] $ is conformally equivalent to a finite cover of the original torus, and hence is also in $\mathcal C _ g $.
\end {proof}

To prove Theorem~\ref {thm:main}, for each positive integer $ g $, one wishes to find spectral data satisfying the conditions of Hitchin's correspondence, (Definition~\ref{spectraldata} and  Periodicity Conditions~\ref{periodicity}).  In particular, the periodicity conditions require  that the periods of the meromorphic differentials $\Theta,\,\Psi $, as well as their integrals over the open curves $ C _ {\pm 1} $
 are all integers.  It suffices to show that the required integrals of $\Theta $ and $\Psi $  each give rational points in $\C P ^{2g+1}$, as then an appropriate multiple of each differential will satisfy the periodicity conditions. In order to guarantee that our tori are rectangular, we will demonstrate the existence of spectral curves
possessing the additional symmetry \mbox{$\lambda\mapsto \frac{1}{\lambda}$}, where \proj.  Using this and other symmetries of the spectral curve, we observe that every spectral curve $\Sigma $ supports differentials $\Theta,\,\Psi $ such that certain of these integrals vanish, whilst the remaining ones are real. Indeed, by assigning to each spectral curve these integrals, one obtains a smooth map $\phi $ from an open subset of the space of spectral curves of genus $ g $ to a product of real projective spaces such that the domain and range of this map have the same dimension.  Our task then reduces to showing that for each $ g $, there is a curve of genus $ g $ at which the differential of this map is invertible.  
An application of the Inverse Function Theorem then yields Theorem~\ref{thm:main}. All linearly full branched minimal immersions $T^2\rightarrow S^3$ are in fact minimal immersions \cite{Hitchin:90}, which is why we are able to conclude that our maps do not have branch points.

\section{Proofs}
As described above, we enforce an additional symmetry \mbox{$\lambda\mapsto \frac{1}{\lambda}$} upon our spectral data in order to find rectangular tori.  This symmetry  induces two
holomorphic involutions on $\Sigma$ and the resulting quotient curves $\Cpm$ are our basic object of study.  We will use proof by
induction, in which at each step the genera of $\Cpm$ increase by one, and hence the genus $g$ of $\Sigma$
increases by two. Thus the proof divides naturally into the even and odd genus cases. We give in this section a  proof for odd genera.  In this case the quotient curves have different genera from one another, which is a complicating factor. The proof for even genera is both simpler and analogous, and is similar to the work of \cite{EKT:93}. A full proof is thus unnecessary, and we give instead a brief description.

\subsection {Odd Genera}
This section forms the technical heart of this manuscript. We begin by describing the algebraic curve data that we shall study, and then turn to the statement of Theorem~\ref{thm:odd}, whose proof is the main purpose of this section.  It states that for each $n\geq 0$, the differential of a certain map $ \phi $ is somewhere invertible. We proceed by induction,
and at each step  pass to the boundary of the moduli space, and then utilise a further limiting argument. We collect in Lemmata~\ref{thm:asym}, \ref{thm:double}, \ref {thm:widehatIpm} the results of the necessary calculations, and derive in Lemma~\ref{thm:extra} an additional assumption that is necessary for the induction step to hold. We then state
Theorem~\ref{thm:oddextra}, which is  simply Theorem~\ref{thm:odd} together with this additional condition, and prove it by induction using  Lemmata~\ref{thm:asym}, \ref{thm:double} and the case $n=0$ (Lemma~\ref {thm:zero}).


Let $Y_+=Y_+(r,x_1,\bar{x_1},\ldots,x_n,\bar{x_n})$ be the curve given by
$$y_{+}^2=(x-r)\prod_{i=1}^{n}(x-x_i)(x-\bar{x_i})$$ and
$Y_-=Y_-(r,x_1,\bar{x_1},\ldots,x_n,\bar{x_n})$ that given by
$$y_{-}^2=(x-2)(x+2)(x-r)\prod_{i=1}^{n}(x-x_i)(x-\bar{x_i})$$
where we assume that \mbox{$r\in (-\infty,-2)\cup(2,\infty)$}, \mbox{$x_i\neq 2$} for \mbox{$i=1,\ldots,2n$} and
\mbox{$x_i\neq x_j$} for $i\neq j$. Denote by \mbox{$\pi_{\pm}:\Cpm\ra\mathbb{C}P^1$} the projections
$(x,y_\pm)\mapsto x$. Construct \mbox{$\pi:\Sigma\ra\mathbb{CP}^1$} as the fibre product of
\mbox{$\pi_{+}:Y_+\ra\mathbb{C}P^1$} and \mbox{$\pi_{-}:Y_-\ra\mathbb{C}P^1$}, that is, let
\mbox{$\Sigma:=\{(p_+,p_-)\in Y_+\times Y_-: \pi_+(p_+)=\pi_-(p_-)\}$} with the obvious projection $\pi$ to $\C
P^1$. Then $\Sigma$ is given by
$$\eta^2=\lambda(\lambda-R)(\lambda-R^{-1})\prod_{i=1}^{n}(\lambda-\lambda_i)(\lambda-{\lambda_i}^{-1})
(\lambda-\bar{\lambda_i})(\lambda-{\bar{\lambda_i}}^{-1}),$$
where $R+R^{-1}=r,\,\alpha_i+{\alpha_i}^{-1}=x_i,$.

$\Sigma$ has genus $2n+1$ and possesses the holomorphic involutions
$$\begin{array}{rccc}
i_{\pm}:&\Sigma&\longrightarrow&\Sigma\\
&(\lambda,\eta)&\longmapsto&\left(\frac{1}{\lambda},\frac{\pm \eta}{\lambda^{2n+1}}\right).
\end{array}$$
The curves $\Cpm$ are the quotients of $\Sigma$ by these involutions, with quotient maps
$$q_{+}(\lambda,\eta)=\Bigl(\lambda+\frac{1}{\lambda},\frac{\eta}{\lambda^{n+1}}\Bigr)=(x,y_+)$$
and
$$q_{-}(\lambda,\eta)=\Bigl(\lambda+\frac{1}{\lambda},\frac{(\lambda+1)(\lambda-1)\eta}{\lambda^{n+2}}\Bigr)=(x,y_-).$$

 $\Cpm$  possesses a real structure $\rho_\pm$, characterised by the
properties that it covers the involution $x\mapsto\bar x$ of $\mathbb{C}P^1$ and fixes the points in
$\pi_\pm\!^{-1}[-2,2]$. These real structures are given by \be \rho_\pm(x,y_\pm)=(\bar x,\mp\bar y_\pm),\mbox{
for $r>2$}\label{eq:r}\ee or
$$\rho_\pm(x,y_\pm)=(\bar x,\pm\bar y_\pm),\mbox{ for $r<-2$.}$$
The cases $r>2$ and $r<-2$ are similar, but the sign difference carries through to future computations. For
simplicity of exposition we assume henceforth that $r>2$. Then the corresponding real structure on $\Sigma$ is
given by
$$\rho(\lambda,\eta)=\Bigl(\frac{1}{\bar \lambda},\frac{-\bar \eta}{{\bar \lambda}^{2n+1}}\Bigr).$$

The curves $Y_\pm$ described above are those which yield spectral curves, but it is easier to work with more general curves that do not possess the real structure $\rho_\pm $ and then identify those in which we are interested by the fact that they satisfy a reality condition.  We thus
 consider also curves \mbox{$\Cpm=\Cpm(r,x_1,\ldots,x_{2n})$} given by

$$y_{+}^2=(x-r)\prod_{i=1}^{2n}(x-x_i)$$ and

$$y_{-}^2=(x-2)(x+2)(x-r)\prod_{i=1}^{2n}(x-x_i)$$ respectively,
where we assume that \mbox{$r\in (2,\infty)$}, \mbox{$x_i\neq\pm 2$} for \mbox{$i=1,\ldots,2n$} and that the
sets \mbox{$\{x_1,x_{2}\},\ldots,\{x_{2n-1},x_{2n}\}$} are mutually disjoint.

\begin{figure}[h]
\hspace*{0cm}\includegraphics{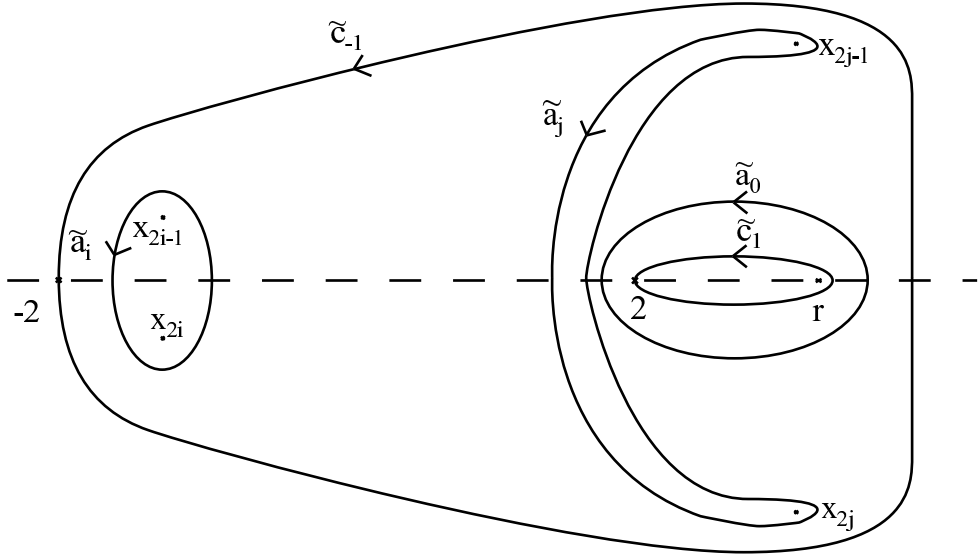} \caption{Curves $\tilde a_i$, $\tilde c_{\pm 1}$ for $p\in
M_{n,\R}$.\label{fig:Cpm}}
\end{figure}

Take $(r,x_1,\ldots,x_{2n})$ as
described above. Let \mbox{$\tilde a_0,\ldots,\tilde a_n $} be simple closed curves in
\mbox{$\cp-\{r,2,-2,x_1,\ldots,x_{2n}\}$}, and $\tilde c_{1}$, $\tilde c_{-1}$ simple closed curves in
\mbox{$\cp-\{r,x_1,\ldots,x_{2n}\}$}, such that (see Figure~\ref{fig:Cpm})
\begin{enumerate}
\item $\tilde a_0$ has winding number one around $2$ and $r$, and winding number zero around the other branch
points of $Y_-$, \item for $i=1\ldots n$, $\tilde a_i$ has winding number one around $x_{2i-1}$ and $ x_{2i}$,
and winding number zero around the other branch points of $Y_\pm$, \item $\tilde c_1$ begins and ends at $x=2$,
has winding number one around $r$ and zero around $x_i$, $i=1,\ldots,2n$, \item $\tilde c_{-1}$ begins and ends
at $x=-2$, and has winding number one around $r$ and each $x_i$ $i=1,\ldots,2n$.
\end{enumerate}
Choose lifts of the curves \mbox{$\tilde a_1,\ldots,\tilde a_n$} to  $Y_+$ and also of \mbox{$\tilde
a_0,\ldots,\tilde a_n$} to $Y_-$. Let \mbox{$a^-_0,a^\pm_1,\ldots,a^\pm_n\in H_1(Y_{\pm},\mathbb{Z})$} denote
the homology classes of these lifts. Denote by \mbox{$b^-_0,b^\pm_1,\ldots,b^\pm_n$} the completions to
canonical bases of  \mbox{$H_1(\Cpm,\mathbb{Z})$.} Choose open curves $c_1$, $c_{-1}$  in $Y_+$ covering the
loops $\tilde c_1$ and $\tilde c_{-1}$.

Denote by $M_n$ the space of $2n+1$-tuples $(r,x_1,\ldots,x_{2n})$ as above
together with the choices we have described. Let $M_{n,\R}$ denote the subset of $M_n$ such that \label{Mn}
\begin{enumerate}
\item $x_{2i}=\bar{x}_{2i-1}$ for $i=1,\ldots,n$, \item for $i=1\ldots n$, $\tilde a_i$ is invariant under
conjugation and intersects the real axis exactly twice, both times in the interval \mbox{$(-2,2)$}, \item the
lifts of  \mbox{$\tilde a_1,\ldots,\tilde a_n$} to $Y_{+}$  are chosen so that the point where $\tilde a_i$
intersects the $x$-axis with positive orientation is lifted to a point where $\frac{y_+}{\im}$ is negative,
\item the lifts of \mbox{$\tilde a_0,\ldots,\tilde a_n$} to $Y_{-}$  are chosen so that the point where $\tilde
a_i$ intersects the $x$-axis with positive orientation is lifted to a point in $Y_-$ where $y_-$ is positive,
\item $c_1$, $c_{-1}$ begin at points with $\frac{y_+}{\im}<0$.
\end{enumerate}

\begin{figure}[h]\hspace*{1.2cm}\includegraphics{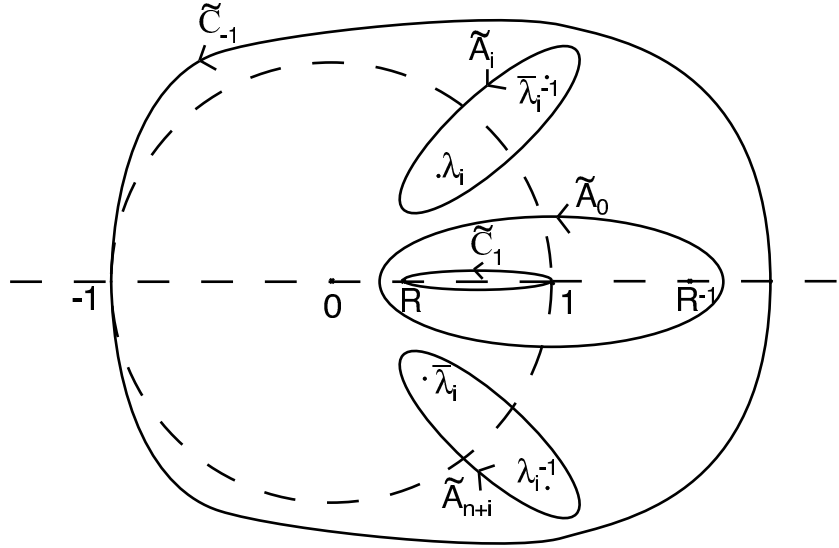}
\caption{The curves $\tilde A_i$ and $\tilde{C}_{\pm 1}$.\label{fig:X}}
\end{figure}

For each $p\in M_{n,\R}$ there is a unique canonical basis \mbox{$A_0,\ldots,A_{2n}$},
\mbox{$B_0,\ldots,B_{2n}$} for the homology of $\Sigma$ such that \mbox{$A_0,\ldots,A_{2n}$} cover the homotopy
classes of loops \mbox{$\tilde A_0,\ldots,\tilde A_{2n}$} shown in Figure~\ref{fig:X} and
$${(q_-)}_*(A_0)=2a^-_0,\;(q_\pm)_*(A_i)=\mp(q_\pm)_*(A_{n+i})=a^\pm_i.$$

There are also unique curves $C_1$ and $C_{-1}$ on $\Sigma$ such that $(q_+)_*(C_{\pm 1})=c_{\pm 1}$; they
project to $\tilde{C}_1$ and $\tilde{C}_{-1}$ of Figure~\ref{fig:X}. Note that $C_1$ connects the two points of
$\Sigma$ with $\lambda=1$ whilst $C_{-1}$ connects the two points of $\Sigma$ with $\lambda=-1$.

Denote by ${\mathcal A}^\pm$ the subgroup of \mbox{$H_1(\Cpm,\mathbb{Z})$} generated by the $a^\pm$ classes.
Then modulo ${\mathcal A}^\pm$,\nopagebreak
$${(q_-)}_*(B_0)\equiv b^-_0,\:(q_\pm)_*(B_i)\equiv
b^\pm_i,\:(q_\pm)_*(B_{n+i})\equiv \mp b^\pm_i.$$
%
%
%
%
%

Let $p\in M_n$ and define $\Wpm=\Wpm(p)$ on $\Cpm(p)$ by:
\begin{enumerate}
\item $\Wpm(p)$ are meromorphic differentials of the second kind: their only singularities are double poles at
$x=\infty$, and they have no residues. \item $\int_{a^-_0}\Omega_-(p) = 0$ and $\int_{a^\pm_i}\Wpm(p) = 0$ for
$i=1,\ldots,n$. \item As $x\ra\infty$, $\Omega_+(p)\ra\frac{x^ndx}{y_+(p)}$ and
$\Omega_-(p)\ra\frac{x^{n+1}dx}{y_-(p)}$.
\end{enumerate}
In view of the defining conditions above, we may write
\be
\Wpm=\frac{\prod_{j= l _\pm}^n(x-\zeta^\pm _j)dx}{y_\pm},\mbox { where } l _ + = 1,\, l _ - = 0.
\label {eq:zeta}
\ee

Let
\be
I_{+}(p):=\im\Bigl(\int_{c_1}\Omega_+(p),\int_{c_{-1}}\Omega_+(p),\int_{b^+_1}\Omega_+(p),
\ldots,\int_{b^+_n}\Omega_+(p)\Bigr),
\label{eq:Ipos}
\ee
\be
I_{-}(p):=\Bigl(\int_{b^-_0}\Omega_-(p),\int_{b^-_1}\Omega_-(p),\ldots,\int_{b^-_n}\Omega_-(p)\Bigr).\label{eq:Ineg}
\ee
Then  for $p\in M_{n,\R}$, $I_+(p)$ and $I_-(p)$ are real, since then
\[
\rho ^*_\pm(\Wpm)=\mp\overline\Wpm
\]
and
$$(\rho_\pm)_*(b^\pm_i)=b^\pm_i \mbox{ mod
${\mathcal A}^\pm$},\:  (\rho_+)_*(c_{\pm 1})=c_{\pm 1} \mbox{ mod ${\mathcal A}^+$},$$ (we have assumed that
(\ref{eq:r}) holds).

Given $p\in M_{n,\R}$, there are real numbers $s_+$ and $s_-$ such that $\im s_+q_+\!^*(\Omega_+(p))$ and
$s_-q_-\!^*(\Omega_-(p))$ are differentials on $\Sigma$ satisfying the conditions of Hitchin's correspondence (Definition~\ref{spectraldata} and Periodicity Conditions~\ref {periodicity}) if
and only if $I_+(p)$ and $I_-(p)$ represent rational elements of $\mathbb{R}P^{n+1}$ and $\mathbb{R}P^{n}$
respectively.

\begin{thm}For each non-negative integer $n$, there exists \mbox{$p\in M_{n,\R}$} such that
\begin{enumerate}
\item $\zeta^+_j(p),\,j=1,\ldots,n$ are pairwise distinct, as are $\zeta^-_j(p),\, j=0,\ldots,n$.
\item The map \vspace*{-5mm}$$\begin{array}{rccc}\\
\phi:&M_n&\longrightarrow&\mathbb{C}P^{n+1}\times\mathbb{C}P^n\\
&p&\longmapsto&([I_+(p)],\,[I_-(p)]).
\end{array}$$
has invertible differential at $p$.
\end{enumerate}
\label{thm:odd}
\end{thm}
This gives that the restriction
$$\phi\left.\right|_{M_{n,\R}}:M_{n,\R}\ra\mathbb{R}P^{n+1}\times\mathbb{R}P^n$$
of $\phi$ to $M_{n,\R}$ also has invertible differential at $p$. The Inverse Function Theorem then implies that
for each positive odd integer $g$, there are countably many spectral curves $\Sigma$ of genus $g$ each giving
rise to a marked torus $(T^2,\tau)$ and a branched minimal immersion \map. The conformal type of the torus is given by (\cite {Hitchin:90})

$$\tau = \frac{\im s_+\rm{p.p.}_\infty(q_+\!^*(\Omega_+(p)))}{s_-\rm{p.p.}_\infty(q_-\!^*(\Omega_-(p)))}=
\frac{\im s_+}{s_-},$$
where $\rm{p.p.}_\infty(q_\pm\!^*(\Omega_\pm))$ denotes the principal part of $q_\pm\!^*(\Wpm)$ at $\infty$.
Thus each torus $(T^2,\tau)$ is rectangular.

In fact we shall prove a slightly stronger result. The extra strength resides in a statement that arises  from
an attempt to prove Theorem~\ref{thm:odd} by induction on $n$, and enables one to complete the induction step.
This statement is somewhat lengthy to formulate, and will appear unmotivated at this juncture. Our approach is
thus to present an attempt to prove Theorem~\ref{thm:odd} by induction, and derive the necessary modifications.
The reader who wishes to view the modified statement at this point is referred to Theorem~\ref{thm:oddextra}, page~\pageref{thm:oddextra}.

\noindent{\bf Proof of Theorem~\ref{thm:odd}:} Suppose then that $\phi$ has invertible differential at
\mbox{$p\in M_{n,\R}$}. For \mbox{$\m\in (-2,2)$}, \mbox{$\n\in\mathbb{R}$}, we shall denote by $(p,\m,\n)$ the
point of \mbox{$M_{n+1,\R}$} such that
$$x_i(p,\m,\n)=x_i(p),\,i=1,\ldots n,\; x_{n+1}(p,\m,\n)=\m+\im\n.$$
Denote by \mbox{$(p,\m,\n)$} the point in \mbox{$M_{n+1,\R}$} with branch points
$$x_i(p,\mu,\nu)=x_i(p),\,i=1,\ldots,2n$$ and
$$x_{2n+1}(p,\mu,\nu)=\mu+\im\nu,\;x_{2n+2}(p,\mu,\nu)=\mu-\im\nu.$$
By considering the boundary case $\n=0$, we show that for a  generic \mbox{$\m\in (-2,2)$} and $\n$ sufficiently
small, $\phi$ has invertible differential at \mbox{$(p,\m,\n)$}.
We shall write \mbox{$\po=\po(p,\m)$} for $(p,\m,0)$.

\begin{figure}[h]
\hspace*{1.5cm}\includegraphics{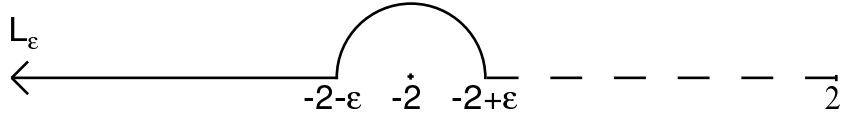} \caption{$h_\epsilon$ is a function of $\m\in
L_\epsilon$.\label{fig:hdomain}}\end{figure}

Let \be H(\m):=\left(\begin{array}{cc}
I_{+}(\po)&0\\
0&I_{-}(\po)\\
\frac{\partial}{\partial r} I_{+}(\po)&\frac{\partial}{\partial r} I_{-}(\po)\\
\frac{\partial}{\partial x_{1}}I_{+}(\po)&\frac{\partial}{\partial x_{1}} I_{-}(\po)\\
\vdots&\vdots\\
\frac{\partial}{\partial x_{2n}}I_{+}(\po)&\frac{\partial}{\partial x_{2n}}I_{-}(\po)\\
\frac{\partial}{\partial \m}
I_{+}(\po)&\frac{\partial}{\partial \m} I_{-}(\po)\\
\frac{\partial^{2}}{\partial {\n}^{2}} I_{+}(\po)&\frac{\partial^{2}}{\partial {\n}^{2}} I_{-}(\po)
\end{array}\right),\label{eq:H}\ee
and \be h(\m):=\det H(\m).\label{eq:he}\ee $h$ is a  real-analytic function of $\m\in(-2,2)$ and for each
\mbox{$\epsilon\in(0,\min_{i=1,\ldots, n}|x_i+2|)$} we may use (\ref{eq:he}) to define it as a real-analytic
function $h_\epsilon$ of $\m$ on the curve $L_\epsilon$ shown in Figure~\ref{fig:hdomain}.

We will show that $h(\m)\not\ra 0$ as $\m\ra\infty$ along  $L_\epsilon$, by computing asymptotics for each of
the vectors in $H(\m)$.
We will also prove that
$$\frac{\partial}{\partial
\n}\left(I_{+}(\po),I_{-}(\po)\right)=0.$$

Then for generic $\mu\in(-2,2)$
 and  $\n$ sufficiently
 small, $d\phi_{(p,\m,\n)}$ is invertible.

A simplification provided by choosing $\n=0$ is that $Y_\pm(p)$ are the respective normalisations of
$Y_\pm(\po)$, with normalisation maps

\begin{center}
$\begin{array}{rccc}
\Psi_\pm:&Y_\pm(p)&\longrightarrow&Y_\pm(\po)\\
&(x,y_\pm(p))&\longmapsto&(x,(x-\m)y_\pm(p))
\end{array}$
\end{center}
and that
%
%
%
%
%
%
%
%

\[
I_\pm (\po)=\Bigl(I_\pm (p),\sqrt {\mp 1}\int_{b^\pm_{n+1}}\Omega_\pm (\po)\Bigr).  \label{eq:Iplus}
\]
%
For each point $q\in M_n$, let $u_\pm(q)$ be local coordinates on $Y_{\pm}(q)$ near $\pi_{\pm}\!^{-1}(\infty)$
such that
$u_{\pm}(q)^{2}=x_{\pm}\!^{-1}$
and
$$y_{\pm}(q)=u_{\pm}(q)^{2 m _\pm +1}+O(u_{+}(q)^{2 m _\pm}) \mbox{ as } x_{\pm}\ra\infty, \mbox { where } m _ + = n,\,m _ -=n +1. $$
Then for $x$ near $\infty$, \be \Wpm(q)=(u_{\pm}(q)+D_{\pm}(q)u_{\pm}(q)^{3}+O(u_{\pm}(q)^{5}))dx,
\label{eq:Dplusdefn}\ee where \be
D_{\pm}(q):=\frac{1}{2}r(q)+\sum_{i=1}^{2n}x_{i}(q)-\sum_{j=l _\pm}^n\zeta^\pm _j(q),\,l _ + =1,\, l _ - = 0.\label{eq:Dplus}\ee

\begin{lem}
As $\m\ra\infty$ along $L_\epsilon$, the following asymptotic expressions hold:
\begin{enumerate}
\item ${\displaystyle I_{\pm}(\po)=\left(I_{\pm}(p),4\sqrt{\mp 1}\m^{1/2}-4\sqrt{\mp 1}
D_{\pm}(p)\m^{-1/2}+O(\m^{-3/2})\right)}$

\item ${\ds\frac{\partial}{\partial{r}} I_\pm (\po)=\biggl(\frac{\partial}{\partial r}
I_{\pm}(p),\sqrt{\mp 1}\Bigl(-2+\sum_{j= l _\pm}^{n}\frac{\partial\zeta^{\pm}_{j}}{\partial
r}\Bigr)\m^{-1/2}+O(\m^{-3/2})\biggr)}$

\item
For $ i=1,\ldots,2n$,\\
${\ds\frac{\partial}{\partial{x_{i}}} I_\pm (\po)=\biggl(\frac{\partial }{\partial x_{i}}
I_{\pm}(p),\sqrt{\mp 1}\Bigl(-2+\sum_{j= l _ \pm}^{n}\frac{\partial\zeta^{\pm}_{j}}{\partial
x_{i}}\Bigr)\m^{-1/2}+O(\m^{-3/2})\biggr)}$

\item ${\ds\frac{\partial }{\partial\m}I_\pm (\po)=\left(0,2\sqrt{\mp 1}\m^{-1/2}+2\sqrt{\mp 1}
D_{\pm}(p)\m^{-3/2}+O(\m^{-5/2})\right)} $

\item ${\ds\frac{\partial}{\partial\n}I_\pm(\po)=0}$
\end{enumerate}
\label{thm:asym}
\end{lem}

\begin {proof} All but the last components of
\setcounter{enumi}{1}({\it\arabic{enumi}})--\setcounter{enumi}{4}({\it\arabic{enumi}}) are applications of
equation~(\ref{eq:Iplus}).  
These last components involve integrals over the curves
$b^{\pm}\!_{n+1}(\po)$. Let $\Gamma$ denote the circle $|x|=\m$, transversed clockwise.  For $\m$ sufficiently
large,
\begin{eqnarray*}\int_{b^+_{n+1}(\po)}\Wp(\po)
&=&-\int_{\Gamma}\frac{\prod_{j=1}^n(x-\zeta^+_j)dx}{\sqrt{(x-r)\prod_{i=1}^n (x-x_i)(x-\bar{x_i})}}
\mbox{ (see Figure~\ref{fig:bnplusone})}\\
&=&4\m^{1/2}-4D_\pm(p)\m^{-1/2}+O(\m^{-3/2}),
\end{eqnarray*}
and similarly for \mbox{$\int_{b^-_{n+1}(\po)}\Omega_-(\po)$}, which gives the remainder of
\setcounter{enumi}{1}({\it\arabic{enumi}})--\setcounter{enumi}{4}({\it\arabic{enumi}}).
\begin{figure}[h]\hspace*{1.5cm}\includegraphics{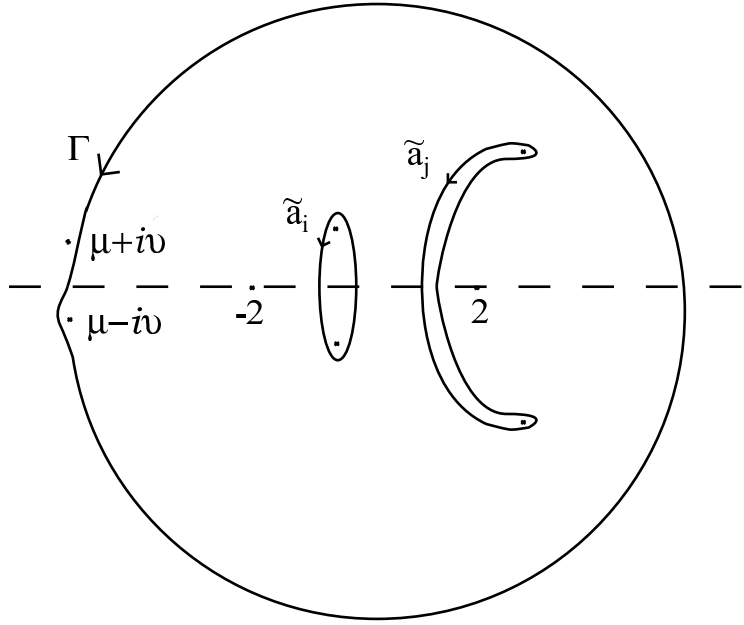}
\caption{We can take a representative of $b^+_{n+1}(\po)$ that projects to a
circle.\label{fig:bnplusone}}\end{figure}

\noindent{\bf\setcounter{enumi}{5}({\it\arabic{enumi}}):} We write $\dot f$ for
$\left.\frac{\partial f}{\partial\n}\right|_{\n=0}$. We will work only with $\Omega_+(\po)$, but similar
arguments apply to $\Omega_-(\po)$. For $i=1,\ldots,n+1$,
\begin{equation}\int_{a^+_i}\dot\Omega_+(\po)=0.\label{eq:dotW}\end{equation}
For $\n$ small, we may write
$$
\Omega_+(p,\m,\n)=\frac{\prod_{j=1}^{n+1}(x-\zeta^+_j(p,\m,\n))dx}{y_+(p,\m,\n)},
$$
 where $\zeta^+_j$ are analytic functions satisfying $$\zeta^+_j(\po)=\left\{\begin{array}{ll}\zeta_j(p)&
\mbox{ for }j=1,\ldots,n\\
\m&\mbox{ for }j=n+1.
\end{array}\right.$$ (Due to this accordance, we write simply $\zeta^\pm_j$ for $\zeta^\pm_j(p)$ or $\zeta^\pm_j(\po)$.)
Then
$$\dot\Omega_+(\po)=\sum_{j=1}^{n+1}\left(\frac{-\dot\zeta^+_j}{x-\zeta^+_j}\right)
\frac{\prod_{k=1}^{n+1}(x-\zeta_k)dx}{y_+}.$$
If $\m$ coincides with one of the $\zeta_j$, $j=1,\ldots,n$, then clearly the lift of
$\dot\Omega_+$ to the normalisation $Y_+(p)$ of $Y_+(\po)$ is holomorphic. If $\m$ does not equal any of the
$\zeta^+_j$, then from (\ref{eq:dotW}), $\dot\Omega_+(\po)$ has zero residue at
$x=\m$, and hence $\dot\zeta_{n+1}=0$ and again $\dot\Omega_+(\po)$ lifts to a holomorphic differential.
This lift has zero $a$-periods,
and so is itself zero. For $i=1,\ldots,n$ then,
$$\left.\pnu\right|_{\n=0}\int_{b^+_i}\Omega_+(\po)=\int_{b^+_i}\dot\Omega_+(\po)=0.$$

To compute $\left.\pnu\right|_{\n=0}\int_{b^+_{n+1}}\Omega_+(\po)$ we use reciprocity with the holomorphic
differential $\omega(p,\m,\n)$ on $Y_+(p,\m,\n)$ defined by
$$\int_{a^+_i}\omega(p,\m,\n) = \Bigg\{\begin{array}{l}0\mbox{, for $i=1,\ldots,n$}\\
2\pi\im\mbox{, for $i=n+1$.}\end{array}$$
To simplify our notation, we write $\ppr$ for  $(p,\m,\n)$.
Reciprocity gives that \be\int_{b^{+}_{n+1}}\Wp(\ppr)=4\ka(\ppr),  \label{eq:intbn+1}\ee
where\be\omega(\ppr)=\frac{\ka(\ppr)\prod_{j=1}^{n}(x-\beta_{j}(\ppr))dx}{y_{+}(\ppr)}.\label{eq:omega}\ee
Since
$$\int_{a^+_i}\dot\omega=\left.\pnu\right|_{\nu=0}\int_{a^+_i}\omega(\ppr)=0,\,i = 1,\ldots, n +1, $$
we have $\dot\omega=0$ . Now $$\dot\omega=\left(\dot
\ka-\ka(\po)\sum_{j=1}^{n}\frac{\dot
\beta_j}{x-\beta_j(\po)}\right)\frac{\prod_{j=1}^{n}(x-\beta^+_j(\po))dx}{y_{+}(\po)},$$ so  $\dot \ka=0$, (and
$ \dot\beta_j=0$, $j=1,\ldots,n$) and thus by (\ref{eq:intbn+1}),
\be\left.\pnu\right|_{\nu=0}\int_{b^{+}_{n+1}}\Wp=0.\label {eq:nuderiv}\ee
\end {proof}
Due to (\ref{eq:nuderiv}), we proceed now to compute $\left.\frac{\partial ^ 2} {\partial\nu ^ 2}\right|_{\nu=0}\int_{b^{+}_{n+1}}\Wp $.

 Using  $\dot\Omega_+(\po)=0$, a
calculation shows that
\[
\ddot{\Omega}_+(\po)=
-\left(\sum_{j=1}^{n}\frac{\ddot{\zeta}^+_j}{x-\zeta^+_j}-
\frac{\ddot{\zeta}^+_{n+1}}{x-\mu}-\frac{1}{(x-\m)^2}\right)\Wp(\po).\label{eq:Wddot}
\]
$\ddot\Omega_+(\po)$ has zero residue at $x=\mu$.
Write $\Wp(\po)=k(x)dx$, then 
$$\ddot\zeta^+_{n+1}k(\mu)=-\frac{dk}{dx}(\mu)$$
and hence using equations~(\ref{eq:Dplusdefn}) and (\ref{eq:Dplus}),%
\begin{equation}\ddot\zeta^+_{n+1}=\frac{1}{2}\mu^ {- 1}+D_+(p)\mu^ {- 2}+O (\mu^ {- 3}).
\label{eq:ddotzeta}\end{equation}
Represent the homology classes $a^+_i,b^+_i$
by loops which project to a
fixed compact region $K\subset\cp $ that is independent of $\m$. For $\mu$ sufficiently large and $x\in K$,
\begin{eqnarray}
(x-\mu)^ {- 1} &=&-\mu^ {- 1}-x\mu^ {-2} +O(\mu^ {- 3})
\label{eq:zminusm}\\
(x-\mu)^ {- 2} &=&\mu^ {- 2}+x\mu^ {-3} +O(\mu^ {- 4})
\end{eqnarray}
and substituting these, together with (\ref{eq:ddotzeta}), into (\ref{eq:Wddot}), we obtain that in $K$,
\be\ddot\Omega_+(\po)=\left(-\sum_{j=1}^n\frac{\ddot\zeta^+_j}{x-\zeta^+_j}-\frac{1}{2\mu^2}
+\frac{2D_+(p)-3x}{2\mu^3}+O\left(\frac{1}{\mu^4}\right)\right)\Wp(\po)\label{eq:ddotW}.\ee
Note that both $\Wp(\po)$ and $\ddot\Omega_+(\po)$ have trivial
$a$-periods, so
$$ \frac{3}{2}\int_{a_i(p)}x\Wp(p)
+\ddot\zeta^+_j\mu^3\sum_{j=1}^n\int_{a^+_i(p)} \frac{\Omega_+(p)}{x-\zeta^+_j}
=O\left(\frac{1}{\mu}\right) \mbox{, for $i=1,\ldots,n$. }$$
Motivated by this, we
define \mbox{$c^\pm _j(p)$} by the equations
\be
\frac{3}{2}\int_{a^\pm _i}x\Wpm(p)+\sum_{j=l _\pm}^{n}c^\pm _j(p)\int_{a^\pm _i}\frac{\Wpm(p)}{x-\zeta^\pm _j}=0,
\:i=l _\pm,\ldots,n,\,l _ + = 1,\, l _ - = 0
\label{eq:cj}\ee
and let
\be
\Wpmh(p):=\frac{3}{2}x\Wpm(p)+\sum_{j= l _\pm}^{n}c^\pm _j(p)\frac{\Wpm(p)}{x-\zeta^\pm _j}
\label{eq:Wph}.\ee
We assumed that
the $\zeta^\pm _j$ are pairwise distinct 
and so the differentials $\frac{\Wpm(p)}{x-\zeta^\pm _j}$ form a
basis for the holomorphic differentials on $Y_\pm (p)$, and $c^\pm_j(p)$ are well-defined.

%
Notice that
\be\ddot\zeta^+_j=c^+_j(p)\mu^ {- 3}+O(\mu^ {- 4}). \label{eq:zetaj}\ee From
 (\ref{eq:ddotW}), and (\ref{eq:zetaj}),
\be\Psi^*_+\ddot\Omega_+= -\frac{1}{2}\mu^ {- 2}\Wp+\mu^ {- 3}\left(D_+(p)\Wp + \Wph\right)
+O\left(\mu^ {- 4}\right).\label {eq:doublederivative}\ee

     %

To compute $\left.\ppnu\right|_{\nu=0}\int_{b^+_{n+1}}\Wp$, we again use reciprocity with $\omega(\po)$.
Differentiating (\ref{eq:omega}), %
one obtains
\be\ddot\omega=\left(\frac{\ddot\ka}{\ka(\po)}-\frac{1}{(x-\mu)^2}
-\sum_{j=1}^{n}\frac{\ddot\beta_j}{x-\beta_j(\po)}\right)\omega(\po).\label{eq:omegaddot}\ee
Taking residues at $x=\mu$, and utilising ${\rm res}_{x=\mu}\omega(\po)=1,\;{\rm res}_{x=\mu}\ddot\omega=0,$ is
gives \be\frac{\ddot\ka}{\ka(\po)}-\sum_{j=1}^{n}\frac{\ddot\beta_j}{x-\beta_j(\po)}={\rm
res}_{x=\mu}\frac{\omega(\po)}{(x-\mu)^2}.\label{eq:ddotka}\ee
Thus in order to obtain an asymptotic expression for $\ddot\ka$, we first obtain expressions for $\ka(\po)$,
$\beta_j(\po)$ and $\ddot\beta_j$. We assume throughout that $x\in K$. We shall
write $\omega(\po)$ for $\Psi_+\!^*(\omega(\po))$.%
 From (\ref{eq:omega}) and (\ref{eq:zminusm}),
%
$$\int_{a^+_i}\frac{\prod_{j=1}^n(x-\beta_j(\po))dx}{y_+(p)}=O\left(\frac{1}{\mu}\right).$$
$\frac{\prod_{j=1}^n(x-\beta_j(\po))dx}{y_+(p)}$ is moreover a differential of the second kind on $Y_+(p)$
whose only singularity is a double pole at $x=\infty$, and it approaches $\frac{x^n}{y_+(p)}$ as $x\ra\infty$. Hence
$$\frac{\prod_{j=1}^n(x-\beta_j(\po))dx}{y_+(p)}+O\left(\frac{1}{\mu}\right)=\Wp(p),$$
\be\beta_j(\po)=\zeta_j+\Oone,\,j=1,\ldots,n.\label{eq:betaj}\ee
%
%
%
%
%
%
%
%
$$\ka(\po)=\mu^{1/2}-D_+(p)\mu^{-1/2}+O(\mu^{-3/2}),$$
and therefore \be\omega(\po)=\mu^{-1/2}\Omega_+(p)+O(\mu^{-3/2}).\label{eq:omega2}\ee
From (\ref{eq:omegaddot}), (\ref{eq:omega2}) and (\ref{eq:zminusm}) we obtain
\[\sum_{j=1}^{n}\ddot\beta_j\int_{a^+_i}\frac{\omega(\po)}{x-\beta_j(\po)}=%
\frac{2}{\mu^{7/2}}\int_{a^+_i}x\Omega_+(\po) + O\left(\frac{1}{\mu^{9/2}}\right) \mbox{ for }i=1,\ldots n.\]
For each $i$, $\int_{a^+_i}x\Omega_+(\po)$ is independent of $\mu$ and so \be
\sum_{j=1}^{n}\ddot\beta_j\mu^{1/2}\int_{a^+_i}\frac{\omega(\po)}{x-\beta_j(\po)}=O\left(\frac{1}{\mu^3}\right).
\label{eq:order3}\ee But by (\ref{eq:betaj}) and (\ref{eq:omega2}),
$$\left(\mu^{1/2}\int_{a^+_i}\frac{\omega(\po)}{x-\beta_j(\po)}\right)^i_j=
\left(\int_{a^+_i}\frac{\Omega_+(p)}{x-\zeta_j}\right)^i_j+O\left(\frac{1}{\mu}\right),$$ and the matrix on the
right is invertible and independent of $\mu$. 
Thus $$\ddot\beta_j(\po) = O\left(\mu^ {- 3}\right),\, j=1,\ldots,n.$$

We now have the asymptotics for $\ka(\po)$, $\beta_j(\po)$ and $\ddot\beta_j$ that we desired earlier;
substituting them into
(\ref{eq:ddotka}) gives
\be\ddot\ka=\frac{3}{8}\mu^{- 3/2}+\frac{9} {8} D_+(p)\mu^{-5/2}+O\left({\mu^{-7/2}}\right).\label {eq:kappaddot}\ee
Set
\be
\widehat{I}_+ (p):=\im\Bigl(\int_{c_1}\widehat{\Omega}_+(p),\int_{c_{-1}}\widehat{\Omega}_+(p),
\int_{b^+_1}\widehat{\Omega}_+(p),\ldots,\int_{b^+_n}\widehat{\Omega}_+(p)\Bigr),\label {eq:Ihatplus}
\ee
\be
\widehat{I}_-(p):=\Bigl(\int_{b^-_0}\widehat{\Omega}_-(p),\int_{b^-_1}\widehat{\Omega}_-(p),
\ldots,\int_{b^-_n}\widehat{\Omega}_-(p)\Bigr).\label {eq:Ihatminus}
\ee
Then utilising (\ref{eq:intbn+1}),  (\ref{eq:kappaddot}) and (\ref{eq:doublederivative}), we have proven the following.
\begin{lem}
${\ds\frac{\partial^{2}}{\partial\n^{2}}I_\pm (\po)=\Bigl(\bigl(-\frac{1}{2}\m^{-2}+D_\pm (p)\m^{-3}\bigr)I_{\pm}(p)-\m^{-3}
\widehat{I}_{\pm}(p)+O(\m^{-4}),}$
\begin{flushright}${\ds\hfill
\frac{3\sqrt{\mp 1}}{2}\m^{-3/2}+\frac{9\sqrt{\mp 1}}{2}D_{\pm}(p)\m^{-5/2}+O(\m^{-7/2})\Bigr)}$\end{flushright}
\label{thm:double}
\end{lem}

We now have asymptotic expressions for each row of the $2n+5\times 2n+5$ matrix $H(\m)$ in (\ref{eq:H}), which
we wish to show is non-singular in the limit as $\m\ra\infty$ along $L_\epsilon$ of Figure~\ref{fig:hdomain}.
The inductive assumption and Lemmata~\ref{thm:asym}, \ref{thm:double} tell us that columns \mbox{$1,\ldots,n+2,n+4,\ldots,2n+4$} of
the the first $2n+3$ rows of $H(\m)$ are linearly independent, and that its $2n+4^{th}$ row $(\pmu I_+(\po);\pmu
I_-(\po))$ is

\noindent${\ds(\underbrace{0,\ldots,0}_{\mbox{\scriptsize $n+2$ zeros}},2\im(\m^{-1/2}+
D_{+}(p)\m^{-3/2})+O(\m^{-5/2});}$
\vspace{-6mm}\begin{flushright}${\ds\underbrace{0,\ldots,0}_{\mbox{\scriptsize $n+1$
zeros}},2(\m^{-1/2}+D_{-}(p)\m^{-3/2})+O(\m^{-5/2})).}$\end{flushright}
Note that the two non-zero entries in this row have leading terms differing only by multiplication by $\im$. We
find a linear combination of the rows of $H(\m)$ that equals\\
${\displaystyle(\underbrace{0,\ldots,0}_{\mbox{\scriptsize $n+2$ zeros}}, \im\m^{-5/2}(4\delta^+(p)-5
D_+(p))+O(\m^{-7/2});}$\vspace{-6mm}\begin{flushright}$\displaystyle{\underbrace{0,\ldots,0}_{\mbox{\scriptsize
$n+1$ zeros}}, \m^{-5/2}(4\delta^-(p)-5D_-(p))+O(\m^{-7/2}))},$\end{flushright} where $\delta^\pm(p)$ are
defined in (\ref{eqn:etas}).
\begin {lem}\label {thm:extra}
 The  matrix $H(\mu)$ is non-singular if
$$\lim_{\m\ra\infty}  4\delta^+(p)-5D_+(p) \neq \lim_{\m\ra\infty}
4\delta^-(p)-5D_-(p)$$ where the limits are taken along $L_\epsilon$.
\end {lem}
\begin {proof} This is the ``extra condition'' referred
to earlier, and we will modify the statement we prove  by  induction to ensure that it is satisfied. First, we
find the linear combination yielding this condition.

From  Lemmata~\ref{thm:asym}, \ref{thm:double} \be\hspace*{-77mm} \Bigl(\ppnu I_+(\po);\ppnu I_-(\po)\Bigr)=\label{eq:ppnu}\ee
$\displaystyle{\Bigl(\bigl(D_+(p)\m^ {- 3}-\frac{1}{2}\m^ {- 2}\bigr)I_+(p) -\m^ {- 3}\widehat
I_+(p)+O(\m ^ {-4}), \frac{3} {2}\im \bigl(\m^{- 3/2}+
3D_{+}\m^{-5/2}+O(\m^{- 7/2})\bigr);}$   \vspace *{-7mm}\begin{flushright}
${\displaystyle (D_-(p)\m^ {- 3}-\frac{1}{2}\m^ {- 2})I_-(p) -\m^ {- 3}\widehat
I_-(p)+O(\m ^ {-4}), \frac{3}{2}\m^{-3/2}+\frac{9} {2}
D_{-}\m^{- 5/2}+O(\m^{- 7/2})\Bigr)}$
\end{flushright}
By the induction hypothesis, there are unique \mbox{$\delta^{\pm}(p)$}, \mbox{$\chi(p)$}, \mbox{$\xi_i(p)$},
\mbox{$i=1,\ldots2n$} such that
\begin{eqnarray}(\widehat{I}_+(p);\widehat I_-(p)) & = & \delta^+(p)(I_+(p);0)+\delta^-(p)(0;I_-(p))
+\chi(p)\frac{\partial}{\partial r}({I}_+(p);I_-(p))\nonumber\\
&&+\sum_{i=1}^{2n}\xi_{i}(p)\pli(I_+(p);I_-(p)).\label{eqn:etas}
\end{eqnarray}
Thus by (\ref{eq:ppnu}) and Lemmata~\ref{thm:asym}, \ref{thm:double} there are
\mbox{$\tilde\delta^{\pm}(p)=\delta^{\pm}(p)+O(\m^{-1})$}, \mbox{$\tilde\chi(p)=\chi(p)+O(\m^{-1})$},
\mbox{$\tilde\xi_i(p)=\xi_i(p)+O(\m^{-1})$} satisfying \vspace*{3mm}

\noindent${\ds\Bigl(\ppnu I_+(\po);\ppnu I_-(\po)\Bigr)
+(\frac{1}{2}\m^ {- 2}+ \bigl(\tilde\delta^+(p)-D_+(p)) \m^ {- 3}\bigr)\left(I_+(\po);0\right)}  $\\*
${\ds +\bigl(\frac{1}{2}\m^ {- 2}+ (\tilde\delta_-(p)-D_+(p))\m^ {- 3}\bigr)\left(0;I_-(\po)\right)
+\tilde\chi(p)\m^ {- 3}\frac{\partial}{\partial
r}(I_+(p);I_-(p))}$\\*
${\ds +\sum_{i=1}^{2n}\tilde\xi_{i}(p)\m^ {- 3}\pli(I_+(p);I_-(p))
=(\underbrace{0,\ldots,0}_{\mbox{\scriptsize $n+2$ zeros}},g_+(\po);
\underbrace{0,\ldots,0}_{\mbox{\scriptsize $n+1$ zeros}},g_-(\po)) \label{eq:lincomb}}$\\*
 where
\[
g_\pm(\m)
=\sqrt {\mp 1}\bigl(\frac{7}{2}\m^{- 3/2}-\frac{3} {2} D_\pm (p)\m^{- 5/2}+4\delta^\pm (p)\m^{- 5/2}
+O(\m^{- 7/2})\bigr)
\]
Denote by $\tilde l(\po)$ the linear combination appearing in (\ref{eq:lincomb}). Define
$ l (\po): =\tilde l (\po) -\frac 74\m^ {- 1}\pmu(I_+(\po),I_-(\po))$; then
\begin{eqnarray}
l (\po) &=& \bigl(\underbrace{0,\ldots,0}_{\mbox{\scriptsize $n+2$ zeros}},
\im(4\delta^+(p)-5 D_+(p))\m^{- 5/2}+O(\m^{- 7/2})\nonumber\\
& &\hspace*{23mm}\underbrace{0,\ldots,0}_{\mbox{\scriptsize $n+1$ zeros}},
(4\delta^-(p)-5D_-(p))\m^{- 5/2}+O(\m^{- 7/2})\bigr)\label{eq:lincomb2}.
\end{eqnarray}
\end {proof}

We are led therefore, to modify the statement we prove by induction to include the assumption that
$$ \lim_{\m\ra\infty}  4\delta^+(p)-5D_+(p) \neq  \lim_{\m\ra\infty}
4\delta^-(p)-5D_-(p),$$%
where the limits are taken along $L_\epsilon$. Of course this modification needs to be such that it is preserved
under the induction step. With this in mind, we define $\delta^\pm(\po)$ by the condition that\\
\be
 (\widehat{I}_+(\po);\widehat{I}_-(\po)) - \delta^+(\po)(I_+(\po);0) +
\delta^-(\po)(0;I_-(\po))\in \label{eq:delta}
\ee
 \vspace{-9mm}\begin{flushright} ${\ds \rm{span}\Bigl\{\pR(I_+(\po);I_-(\po)),
\pli (I_+(\po);I_-(\po)),\pmu (I_+(\po);I_-(\po)), \ppnu (I_+(\po);I_-(\po))\Bigr\}}$\end{flushright}
\vspace{0mm} and calculate the relationship between $\delta^+(\po)-\delta^-(\po)$ and $\delta^+(p)-\delta^-(p)$.

\begin{lem}
As $\m\ra\infty$ along $L_\epsilon$,
\[
\widehat I_\pm (\po)=(\widehat I_\pm (p),\sqrt {\mp 1} (2\m^{3/2}+
(6 D_\pm (p)-4\delta^\pm (p))\m^{1/2} + O(\m^{-1/2}))
\]
\label{thm:widehatIpm}\end{lem}
\begin {proof} For $\m$ sufficiently large, we may assume that
$\m\neq\zeta^\pm_j$,\\$j=1,\ldots,n$. Then (again arguments for $\Wm$ are similar to those for $\Wp$)
$$\widehat\Omega_+(\po):=\frac{3}{2}x\Wp(\po)
+\sum_{j=1}^{n+1}\frac{c^+_j(\po) \Wp(\po)}{x-\zeta^+_j},$$
where the $c^+_j(\po)$, $j=1,\ldots,n+1$ are determined by the (non-singular)
system of equations
$$\frac{3}{2}\int_{a^\pm_i}x\Wpm(\po)+\sum_{j=1}^{n+1}c^\pm_j(\po)\int_{a^\pm_i}\frac{\Wpm(\po)}{x-\zeta^\pm_j}
=0,\:i,\,j=1,\ldots,n+1.$$
Taking $i=n+1$ we  see that \mbox{$c^+_{n+1}(\po)=0$} and \mbox{$c^+_j(\po)=c^+_j(p)$,} for
\mbox{$j=1,\ldots,n$}, so
$$\Psi_+\!^*(\widehat \Omega_+(\po))=\widehat \Omega_+(p).$$
proving all but the last component of the
lemma above. For this, we again let  $\Gamma$ denote the circle $|x|=\m$, traversed clockwise. For $\m$
sufficiently large,
\begin{eqnarray*}\int_{b^+_{n+1}}\widehat\Omega_+(\po)
&=&-\int_{\Gamma}\biggl(\frac{3}{2}x+\sum_{j=1}^{n}\frac{c^+_j(p)}{x-\zeta^+_j}\biggr)
\frac{\prod_{k=1}^n(x-\zeta^+_k)dx}{\sqrt{(x-r)\prod_{i=1}^{2n} (x-x_i)}}\\
&=&2\im\m^{3/2}+(6\im D_+(p)-4\im\delta^+(p))\m^{1/2}+O(\m^{-1/2}),
\end{eqnarray*}
and similarly for \mbox{$\int_{b^-_{n+1}(\po)}\widehat\Omega_-(\po)$}.
\end {proof}
From Lemmata~\ref{thm:asym}, \ref{thm:double}, \ref{thm:widehatIpm},

\noindent${\ds (\widehat I_+(\po);\widehat I_-(\po))-\delta^+(p)(I_+(\po);0)
-\delta^-(p)(0;I_-(\po))-\chi(p)\pR(I_+(\po);I_-(\po))}$\\*
${\ds
-\sum_{i=1}^{2n}\xi_i(p)\pli(I_+(\po);I_-(\po))}$\\* \hspace*{5mm}${\ds=(0,2\im\m^{3/2}+(6\im
D_+(p)-4\im\delta^+(p))\m^{1/2}+O(\m^{-1/2});}$ \vspace{-0mm}\begin{flushright}${\ds 0,
2\m^{3/2}+(6D_-(p)-4\delta^-(p))\m^{1/2}+O(\m^{-1/2})}$\end{flushright}
\vspace{-0mm} \be\hspace*{-59mm}=\Lambda \pmu(I_+(\po);I_-(\po)) + \Upsilon l(\po),\label{eq:LamUps}\ee
where $l(\po)$ is defined in (\ref{eq:lincomb2}) and $\Lambda$ and
$\Upsilon$ are defined by the equations
\begin{flushleft}${\ds \left(\begin{array}{cc}
\frac{2\im}{\m^{1/2}}(1 + \frac{D_+(p)}{\m})+O(\frac{1}{\m^{5/2}}) &
\frac{\im(4\delta^+(p)-5D_+(p))}{\m^{5/2}}+O(\frac{1}{\m^{7/2}})\\
\frac{2}{\m^{1/2}}(1 + \frac{D_-(p)}{\m})+O(\frac{1}{\m^{5/2}}) &
\frac{4\delta^-(p)-5D_-(p)}{\m^{5/2}}+O(\frac{1}{\m^{7/2}})
\end{array}\right)
\left(\begin{array}{c}
\Lambda\\
\Upsilon
\end{array}\right)}$\end{flushleft}\nopagebreak
\begin{flushright}${\ds =\left(\begin{array}{c}
2\im\m^{3/2}+ (6\im D_+(p)-4\im\delta^+(p))\m^{1/2} + O(\frac{1}{\m^{1/2}})\\
2\m^{3/2}+ (6 D_-(p)-4\delta^-(p))\m^{1/2} + O(\frac{1}{\m^{1/2}}).
\end{array}\right)}$\end{flushright}
Hence
using (\ref{eq:lincomb2}) and (\ref{eq:LamUps}),
\[\delta^+(\po)-\delta^-(\po)=
\frac{(D_+(p)-D_-(p))\left(3(\delta^+(p)-\delta_-(p))-4(D_+(p)-D_-(p))\right)}
{4(\delta^+(p)-\delta_-(p))-5(D_+(p)-D_-(p))}.\label{eq:etapo} \]

Defining $T_p$ to be the linear fractional transformation
$$T_p:x\mapsto (D_+(p)-D_-(p))\frac{3x-4(D_+(p)-D_-(p))}{4x-5(D_+(p)-D_-(p))},$$
$$T_p(\delta^+(p)-\delta^-(p))=\delta^+(\po)-\delta^-(\po).$$
Moreover, we know that $D_\pm(\po)=D_\pm(p)$ and thus
$$T_{\po}=T_p.$$
In order to conclude the proof of Theorem~\ref{thm:odd},  it suffices then to show:
\begin{thm}
    For each positive integer $m$ and integer $n$ with $0\leq n\leq m$ there exists $p\in M_{n,\R}$ such that
\begin{enumerate}
\item{$\zeta^+_j(p)$, $j=1,\ldots,n$ are pairwise distinct, as are $\zeta^-_j(p)$, $j=0\ldots n$,} \item
$\mathbb R^{2n+3}$ is spanned by the vectors $(I_{+}(p),0)$, $(0,I_{-}(p))$, $\frac{\partial}{\partial
r}\left(I_{+}(p), I_{-}(p)\right)$
 and $\pli\left(I_{+}(p), I_{-}(p)\right)$, $i=1,\ldots,2n$,
\item{$5\left(D_{+}(p)-D_{-}(p)\right)+4T^{k}_p\left(\delta_{+}(p)-\delta_{-}(p)\right)\neq 0$, for $0\leq k\leq
m-n$.}
\end{enumerate}
\label{thm:oddextra}
\end{thm}%
$M _ {n,\R} $ is defined on page~\pageref{Mn},  $\zeta ^\pm,\,I _+,\,I_-,\,D _\pm $ are defined in equations (\ref{eq:zeta}), (\ref{eq:Ipos}), (\ref{eq:Ineg}), (\ref{eq:Dplus}) respectively and $\delta ^\pm $ is given by equations (\ref{eq:cj}), (\ref{eq:Wph}), (\ref{eq:Ihatplus}), (\ref{eq:Ihatminus}) and (\ref{eq:delta}).
\begin {proof} Fix $m$, and for $n<m$ suppose $p\in M_{n,\R}$ satisfies the
conditions of Theorem~\ref{thm:oddextra}. By the above arguments, the set of $\m\in(-2,2)$ such that
\begin{enumerate}
\item[\setcounter{enumi}{1}({\roman{enumi}})] for all $\epsilon\in (0,\min_{i = 1,\ldots,n}|x_i+2|)$, \mbox{$
h_\epsilon(\m)\neq 0$} (see (\ref{eq:he})),

\item[\setcounter{enumi}{2}({\roman{enumi}})] for $j=1,\ldots n$, $\m\neq\zeta^+_j$ and

\item[ \setcounter{enumi}{3}({\roman{enumi}})]for $j=0,\ldots n$, $\m\neq\zeta^-_j$
\end{enumerate}
is dense  in $(-2,2)$.

Take such a $\m$. Then $\po=(p,\m,0)$ satisfies Theorem~\ref{thm:oddextra}, where in
\setcounter{enumi}{2}({\it\arabic{enumi}}) we replace \mbox{$\frac{\partial}{\partial x_{2n+2}}(I_+(\po);
I_-(\po))$} by \mbox{$\ppnu(I_+(\po); I_-(\po))$}. Then for $\n$ small,
$$\delta^\pm(p,\m,\n) = \delta^\pm(\po) + O(\n),$$
$$D_\pm(p,\m,\n) = D_\pm(\po) +O(\n)$$
and
$$T_{(p,\m,\n)}=T_{\po}+O(\n),$$
we conclude that $(p,\m,\n)$ satisfies Theorem~\ref{thm:oddextra}. It remains to show the existence of $p\in
M_{0,\R}$ verifying \setcounter{enumi}{1}({\it\arabic{enumi}}) and \setcounter{enumi}{2}({\it\arabic{enumi}}) of
Theorem~\ref{thm:oddextra}, and such that for no $k\geq 0$ do we have
$5\left(D_{+}(p)-D_{-}(p)\right)+4T^{k}_p\left(\delta_{+}(p)-\delta_{-}(p)\right)=0$.

\noindent{\bf Genus One ($n=0$)}
We consider pairs $Y_+=Y_+(r)$ and $Y_-=Y_-(r)$ given by
$$y_+^2=(x-r)$$
and
$$y_-^2=(x+2)(x-2)(x-r)$$ respectively, where $r>2$. Writing
$\pi_\pm:(x,y_\pm)\mapsto x$ for the projections to $\cp$, the fibre product of these is the genus one curve
$\Sigma=\Sigma(r)$, given by
$$\eta^2=\lambda(\lambda-R)(\lambda-\frac{1}{R})\mbox{, where } R+\frac{1}{R}=r.$$
\begin{lem}\label{thm:zero}
There exists a $p\in M_{0,\R}$ such that
\begin{enumerate}
\item $\mathbb{R}^{3}$ is spanned by the vectors $(I_+(p),0),\,(0,I_-(p))$ and $\pR(I_+(p),I_-(p))$, \item for
all $k\geq 0$, $5(D_+(p)-D_-(p))+4T^k(p)(\delta^+(p)-\delta^-(p))\neq 0$.
\end{enumerate}
\label{lem:one}
\end{lem}
\begin {proof} The natural limit to consider is $R\ra 1$, i.e. $r=R+1/R \ra 2$,
 which suggests setting $\zeta:=x+2$, $t:=r-2$. Then $Y_+(t)$ is given by $$y_+^2:=\z-t,$$
and $Y_-(t)$ by
$$y_-^2:=\z(\z-t)(\z+4).$$
For each $t>0$, choose $c_1(t), c_{-1}(t)$ and $a^-(t)$ as shown in Figure~\ref{fig:oneCpm}.
\begin{figure}[h]
\hspace*{1.5cm}\includegraphics{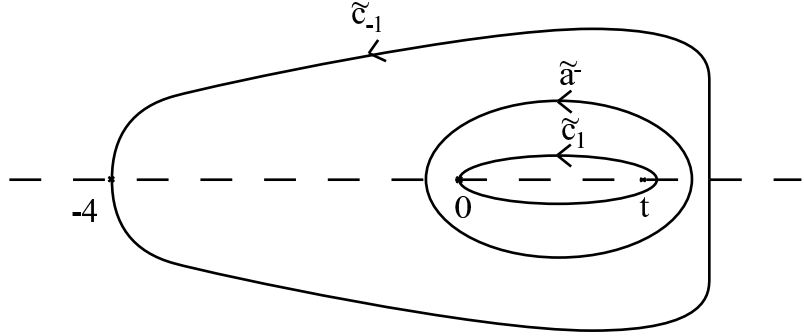} \caption{ The curves $\tilde a^-$ and $\tilde c_{\pm 1}$, for
$n=0$.\label{fig:oneCpm}}
\end{figure}
We write $\Wm(t)=\frac{(\z-s(t))d\z}{y_-}$, where $s(t)$ is defined by the condition $\int_{a^-(t)}\Wm(t) = 0$.
Since
$$\int_{a^-(t)}\Wm(t) = 0\;\mbox{ and }\; \frac{\partial}{\partial t}\int_{a^-(t)}\Wm(t) = 0, $$
we have \be s(t)=t+O(t^2),\label{eq:s}\ee and so
\[ I_-(t)%
=8+O(t).\]
Now $I_+(t)=\im(\int_{c_1(t)}\Wp(t),\,\int_{c_{-1}(t)}\Wp(t))$, where $c_1(t)$ is a path in $Y_+(t)$ joining the
two points with $\z=0$, and $c_{-1}(t)$ is one joining the two points with $\z=-4$, both beginning at points
with $\frac{y_+}{\im}<0$. Then
\[ I_+(t)%
=(4t^{1/2},4(4+t)^{1/2}) \] and
$$\frac{\partial I_+(t)}{\partial t} = (2t^{-1/2},2(4+t)^{-1/2})$$ so we see that condition (1)
of Lemma~\ref{lem:one} is satisfied for all $r>2$.

 We have
$\lim_{t\ra 0}D_+(t)%
=1$, $\lim_{t\ra 0}D_-(t)
=-1$ and we proceed  to calculate  $\lim_{t\ra 0}\delta^\pm(t)$.
$$\widehat\Omega_+(t)=\frac32x\Wp(t) = \frac{3(\z+2)d\z}{2\sqrt{\z-t}},$$ so
\[
\widehat I_+(t)%
=(12t^{1/2}+8t^{3/2},\,12(4+t)^{1/2}+8(4+t)^{3/2}).  \] Recall that
$\widehat\Omega_-(t)=\frac{3}{2}(\z+2)\Wm(t)+\frac{c_-(t)\Wm(t)}{\z-s(t))}$, where $c_-(t)$ is defined by
\begin{equation}\frac{3}{2}\int_{a^-(t)}\frac{(\z+2)(\z-s(t))d\z}{\sqrt{\z(\z+4)(\z-t)}}+
c_-(t)\int_{a^-(t)}\frac{d\z}{\sqrt{\z(\z+4)(\z-t)}}=0. \label{eq:cminus}\end{equation}
Residue calculations show that $$\int_{a^-(t)}\frac{d\z}{\sqrt{\z(\z+4)(\z-t)}}=\pi\im -\frac{\pi\im}{16}t +
O(t^2),$$
$$\int_{a^-(t)}\frac{\z d\z}{\sqrt{\z(\z+4)(\z-t)}}=\pi\im t +O(t^2)$$ and
$$\int_{a^-(t)}\frac{\z^2 d\z}{\sqrt{\z(\z+4)(\z-t)}}=O(t^2).$$
Substituting these and (\ref{eq:s}) into (\ref{eq:cminus}), gives
$$c_-(t)=O(t^2).$$
so, using (\ref{eq:s}), we find that \[
\widehat I_-(t)%
=-8+O(t). \] Thus
\begin{eqnarray*}
\left(\begin{array}{cc}
I_+(t)&0\\
0&I_-(t)\\
\frac{\partial I_+(t)}{\partial t}&\frac{\partial I_-(t)}{\partial t}\\
\widehat I_+(t)&\widehat I_-(t)
\end{array}\right)&=&\left(\begin{array}{ccc}
4t^{1/2}&4(4+t)^{1//2}&0\\
0&0&8+O(t)\\
2t^{-1/2}&2(4+t)^{-1/2}&O(1)\\
12t^{1/2}+8t^{3/2}&12(4+t)^{1/2}+8(4+t)^{3/2}&-8+O(t)
\end{array}\right).
\end{eqnarray*}
Upon multiplication of its third row by $t$, its first column by $2t^{-1/2}$ and its second column by
\(2(4+t)^{-1/2}\) this matrix becomes
\[\left(\begin{array}{ccc}
2&2&0\\
0&0&8\\
1&0&0\\
6&22&-8
\end{array}\right)+O(t).\]
Since \(-11(2,2,0)+1(0,0,8)+16(1,0,0)+1(6,22,-8)=(0,0,0)\), then recalling that \(\delta^\pm(t)\) are defined by
the condition
$$(\widehat I_+(t),\widehat I_-(t))+\delta^+(t)(I_+(t),0)+\delta^-(t)(0,I_-(t))\in {\rm span}
\left\{\frac{\partial}{\partial t}(I_+(t),I_-(t))\right\},$$
we conclude that
$$\lim_{t\ra 0}\delta^+(t)=-11$$ and
$$\lim_{t\ra 0}\delta^-(t)=1.$$
The linear fractional transformation $T_t$ is defined by
$$T_t:u\mapsto \frac{-(D_-(t)-D_+(t))(3u+4(D_-(t)-D_+(t)))}{4u+5(D_-(t)-D_+(t))},$$
so letting $T:=\lim_{t\ra 0}T_t$,
$$T:u\ra \frac{3u-8}{2u-5}.$$
This has a  unique fixed point ($u=2$) and so is conjugate to a translation, in fact denoting the map
$u\mapsto\frac{1}{u-2}$ by $S$, we have
$$STS^{-1}:u\mapsto u-2.$$
Now
\[ \begin{array}{lccc}
&4T^k(\lim_{t\ra 0}(\delta^-(t)-\delta^+(t)))&=&5(\lim_{t\ra 0}(D_+(t)-D_-(t)))\\
\Leftrightarrow &T^k(12)&=&\frac52\\
\Leftrightarrow &(STS^{-1})^k(\frac{1}{10})&=&2\\
\Leftrightarrow &\frac{1}{10}-2k&=&2,
\end{array}\]
which is clearly false for all integers $k\geq 0$. Thus for $t>0$ sufficiently small, the Lemma holds.
\end {proof} This concludes the proof of Theorem~\ref {thm:oddextra}, and hence also that of Theorem~\ref{thm:odd}.\end {proof}

\subsection {Even Genera}
We give a brief indication of how to prove Theorem~\ref{thm:main} for even genera.  In this case the quotient curves have the same genus, and the proof is both simpler than that for the odd genus case, and similar to the proof appearing in \cite{EKT:93}. For these reasons, we content ourselves with describing the appropriate even genus analogue of Theorem~\ref{thm:odd}, whose statement and proof was the main purpose of the previous section.

One is now interested in spectral curves $\Sigma $ of the form
$$\eta^2=\lambda\prod_{i=1}^{n}(\lambda-\lambda_i)
(\lambda-{\lambda_i}^{-1})(\lambda-\bar{\lambda_i})(\lambda-{\bar{\lambda_i}}^{-1}).
$$
 These possess  a real structure
$$\rho:(\lambda,\eta)\mapsto({\bar \lambda^{-1}},{\bar \eta}{\bar \lambda^{-(2n+1)}})$$
and holomorphic involutions
$$\begin{array}{rccc}
i_{\pm}:&\Sigma&\longrightarrow&\Sigma\\
&(\lambda,\eta)&\longmapsto&\left(\frac{1}{\lambda},\frac{\pm \eta}{\lambda^{2n+1}}\right).
\end{array}$$
The  quotients of $\Sigma$ by these involutions have the same genus, and are given by
$$y_\pm^2=(x\pm 2)\prod_{i=1}^{2n}(x-x_i), $$
 with quotient maps
$$q_\pm(\lambda,\eta)=\left(\lambda+\frac{1}{\lambda},\frac{(\lambda\pm 1)\eta}{\lambda^{n+1}}\right)=(x,y_\pm).  $$
The real structure on $\Sigma $ induces real structures
$$\rho_\pm(x,y_\pm)= (\bar x,\pm\bar y_\pm)$$ on $\Cpm$.
One can define a standard homology basis $A_1,\ldots,A_{2n}$,
 $ B_1,\ldots,B_{2n}$ for $\Sigma $ and open curves $ C _ {1},\, C _ {-1}  $, where $C_{\pm 1} $ joins the two points on $\Sigma $ with $\lambda =\pm 1 $, such that
\begin {enumerate}
\item $(q_\pm)_*(A_i)=\pm(q_\pm)_*(A_{n+i}) $, $i = 1,\ldots, n $
\item $ \rho _*(B _i)\cong B_i\mbox{ mod } \mathcal A,\,\rho _*(C_{\pm})\cong C_{\pm}\mbox{ mod } \mathcal A $,
\end {enumerate}
where $\mathcal A$ is the subspace of $H_1(\Sigma,\Z)$ generated by the $A_j$.
\begin{figure}[h]
\hspace*{1.5cm}\includegraphics{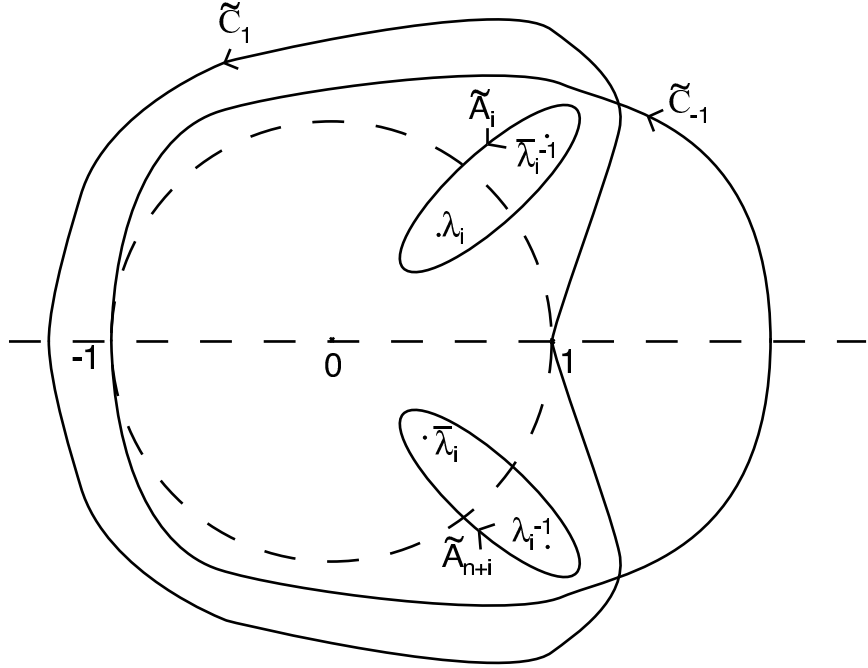} \caption{Projections $\tilde A_i$ and $\tilde {C}_{\pm 1}$ of
$A_i$ and $ {C_{pm 1}}$ to the $\lambda $-plane.\label{fig:evenX}}
\end{figure}
Set
\[
a ^\pm_i=(q_\pm)_*(A_i),\,b^\pm_i=(q_\pm)_*(B_i)
\]
and define differentials $\Wpm=\Wpm(p)$ on $\Cpm(p)$ by:
\begin{enumerate}
\item $\Wpm(p)$ are meromorphic differentials of the second kind: their only singularities are double poles at
$x=\infty$, and they have no residues.
\item $\int_{a^\pm_i}\Wpm(p) = 0$ for $i=1,\ldots,n$. \item As
$x\ra\infty$, $\Wpm(p)\ra\frac{x^n dx}{y_\pm(p)}$.
\end{enumerate}
Then
\[
\rho ^*_\pm (\Wpm)=\pm\overline\Wpm
\]
so
$$I_{\pm}(p):=\sqrt{\pm 1}\Bigl(\int_{c_ {\pm 1}}\Omega_\pm(p),\int_{b^\pm _1}\Omega_\pm (p),\ldots,\int_{b^\pm _n}\Omega_\pm(p)\Bigr)$$
 are real.  If they give rational elements of $\R P^{n}$, then certain real multiples of $\sqrt{\pm 1}(q_\pm) ^*(\Omega _\pm)$ will satisfy the periodicity conditions (\ref{periodicity}).

Thus the natural analog of Theorem~\ref{thm:odd} is the statement that for each $n$, there is a point at which the map
\[
(x_1,\ldots, x_{2n})\mapsto ([I_+], [I_-])
\]
has invertible differential.  This can again be proven by induction, although again additional conditions are required to yield the induction step.

\bibliographystyle{plain}
\bibliography{harmonic}
\end{document}

%% file: defns.tex
\newcommand{\map}{\mbox{$f:(T^2,\tau)\rightarrow S^3$ }}
\newcommand{\proj}{\mbox{$\pi:\Sigma\rightarrow \mathbb{C}P^1$}, \mbox{$(\lambda, \eta)\mapsto \lambda$ }}
\newcommand{\ra}{\rightarrow}

\newcommand{\Wpm}{\Omega_{\pm}}
\newcommand{\Wp}{\Omega_+}
\newcommand{\Wm}{\Omega_-}
\newcommand{\Cpm}{Y_{\pm}}

\newcommand{\cp}{\mbox{$\mathbb{C}P^1$}}
\newcommand{\Wpmh}{\mbox{$\widehat{\Omega}_{\pm}$ }}
\newcommand{\Wph}{\widehat\Omega_+}

\newcommand{\pR}{\frac{\partial}{\partial r}}
\newcommand{\pmu}{\frac{\partial}{\partial\mu}}
\newcommand{\pnu}{\frac{\partial}{\partial\nu}}
\newcommand{\ppnu}{\frac{\partial^2}{\partial\nu^2}}
\newcommand{\pli}{\frac{\partial}{\partial x_i}}

\newcommand{\m}{\mu}
\newcommand{\n}{\nu}
\newcommand{\im}{\sqrt{-1}}
\newcommand{\be}{\begin{eqnarray}}
\newcommand{\ee}{\end{eqnarray}}
\newcommand{\bes}{\begin{eqnarray*}}
\newcommand{\ees}{\end{eqnarray*}}
\newcommand{\p}{p^\prime}
\newcommand{\po}{\p_0}
\newcommand{\ka}{\kappa}

\newcommand{\Oone}{O\left(\frac{1}{\mu}\right)}
\newcommand{\R}{\mathbb{R}}
\newcommand{\ds}{\displaystyle}
\newcommand{\z}{\zeta}

\newcommand{\ppr}{p^\prime}